\documentclass[11pt,psfig]{article}
\usepackage{amssymb,amsmath,amscd}
\usepackage[all]{xy}
\usepackage{mathtools}
\newtheorem{lem}{Lemma}[section]
\newtheorem{thm}{Theorem}[section]
\newtheorem{cor}{Corollary}[section]

\newcommand{\mb}{\mathbb}

\newcommand{\commentout}[1]{}

\newcommand{\mc}{\mathcal}
\newcommand{\ep}{\epsilon}
\newcommand{\arr}[1]{\left( \begin{array}{clcr} #1 \end{array} \right)}

\newcommand{\sgn}{\, {\rm sgn}}

\DeclarePairedDelimiter\ceil{\lceil}{\rceil}

\newcommand{\vertiii}[1]{{\left\vert\kern-0.25ex\left\vert\kern-0.25ex\left\vert #1 
    \right\vert\kern-0.25ex\right\vert\kern-0.25ex\right\vert}}

\begin{document}
\title{ A Thin Fundamental Set for $SL(2, \mathbb Z)$}
\author{Hongyu He \footnote{Key word:  Iwasawa decomposition,  Siegel set, $SL(2, \mathbb Z)$, fundamental domain} \\
Department of Mathematics \\
Louisiana State University \\
email: hhe@lsu.edu\\
}
\date{}
\maketitle

\abstract{Let $\Gamma=SL(2, \mathbb Z)$ and $G=SL(2, \mathbb R)$.   Let $g=kan$ be the Iwasawa decomposition. Let $\ep$ be a small positive number.  In this paper, we construct a fundamental set $\mc F_{\epsilon}$ such that the $k$-component of $ g \in \mc F_{\epsilon}$ is  within the $\epsilon$-distance from the identity. We further prove an inequality for the 
$L^2$-norm of functions on $G/\Gamma$.}

\section{Introduction}
We start with the projective group $PSL(2, \mathbb R)$. Any element $g \in PSL(2, \mb R)$  has an Iwasawa decomposition $kan$ with $k \in PSO(2), a \in A, n \in N$, where $A$ consists of diagonal matrices with positive entries $(a, a^{-1})$ and $N$ consists of upper triangular unipotent matrices parametrized by $t \in \mathbb R$. Let  $PSL(2, \mathbb Z)$ be the modular group consisting of all matrices in $PSL(2, \mathbb R)$ with integer entries. Automorphic forms on $PSL(2, \mathbb R)/PSL(2, \mathbb Z)$ play a central role in many branches of mathematics \cite{borel}. Their analytic properties were often obtained by analysis on the fundamental set $\mc F$, explicitly
$$\{ kan: k \in PSO(2), |t| \leq \frac{1}{2}, a^{-4}+t^2 \geq 1. \}$$
Here the fundamental set has a cusp at $0$. This is consistent with \cite{hc} and \cite{borel19}, but differs from the one more commonly used by an inversion $a \rightarrow a^{-1}$ (\cite{go}). \\
\\
One advantage of using $\mc F$ as the fundamental set, is that analysis based on $\mc F$ will often involve computations on $K$-finite functions which can be expressed as hypergeometric functions ${}_{2}F_{1}$ or ${}_1 F_{1}$ (\cite{val}). However, the behavior of hypergeometric functions can be very complicated and precise computations are often impossible (\cite{ba}). 
In this paper, we shall construct a fundamental set that is not $K$-invariant. This fundamental set $\mc F_{\epsilon}$ will only involve a small neighborhood of the compact group $K$. This small neighborhood of $K$ can be made infinitesimally small.\\
\\
To state our main result, we let $G=SL(2, \mathbb R)$ and $\Gamma=SL(2, \mathbb Z)$. Fix the standard Iwasawa decomposition $KAN$ with $N$ the unipotent upper triangular matrices parametrized by $t \in \mathbb R$, $K = SO(2)$ parametrized by $\theta \in \mathbb R/\mathbb Z$. Our main result can be stated as follows.
\begin{thm} Let $\mc F_{\epsilon}= \cup_{i=1}^4 \mc F^i_{\epsilon} \subseteq G$ with
$$\mc F^1_{\epsilon}= \{ g=kan: |\theta| < \epsilon, a^2 \leq \csc(\epsilon+|\theta|), |t| \leq \frac{1}{2} \};$$
$$\mc F^2_{\epsilon}=\{ g=kan: |\theta| < \epsilon, \csc(\epsilon+|\theta|) \leq a^2 \leq \csc(\epsilon-|\theta|), \sgn(\theta) t \in [-a^{-2}\cot(\epsilon+|\theta|),(1-\sqrt{1-a^{-4}})] \};$$
$$ \mc F^3_{\epsilon}=\{ g=kan: |\theta| < \epsilon, \csc(\epsilon-|\theta|) \leq a^2 \leq \cot(\epsilon-\theta)+\cot(\epsilon+\theta), t \in [1-a^{-2} \cot(\epsilon+\theta), a^{-2} \cot(\epsilon-\theta)] \}; $$
$$ \mc F^4_{\epsilon}=\{ g=kan: |\theta| < \epsilon, \csc(\epsilon-|\theta|) \leq a^2 , t \in [\sqrt{1-a^{-4}}-1, 1-\sqrt{1-a^{-4}}] \}.$$
The natural map $\pi|_{\mc F_{\ep}}: \mc F_{\epsilon} \subseteq G \rightarrow G/\Gamma$ is surjective. $\pi|_{\mc F_{\ep}}$ is injective on the interior of $\mc F_{\ep}$ and finite on the boundary of $\mc F_{\ep}$.
\end{thm}
The fundamental set $\mc F_{\epsilon}$ is contained entirely in $K_{\ep} AN$ with $K_{\epsilon}$ the segment of $SO(2)$ within $\epsilon$-distance from the identity.  See Theorem \ref{big} for a stronger statement.\\
\\
The main idea of the proof is to interpret $G/\Gamma$ as the space of unital lattices in $\mathbb R^2$, namely those lattices with two generators that span an area of $1$. Every fundamental set for $G/\Gamma$ essentially corresponds to a parametrization of the unital lattices. In the classical reduction theory, the parametrization of a unital lattice is based  on the minimal element in the unital lattice (\cite{borel19}). In this paper,  we modify the classical reduction theory by imposing a restriction on the $K$-component. Then the relations between $a$-component and $t$-component become a lot more complicated, but still tractable. We prove a sufficient and necessary condition that dictates the relation between $a$ and $t$. This leads us to the construction of $\mc F_{\epsilon}$. \\
\\
At the end of this paper, we also prove some integral inequality similar to \cite{he2}. These bounds relate the $L^2$-norm on certain \lq\lq conic\rq\rq region to the $L^2$-norm on the fundamental set. \\
\\
Finally, we shall remark that most results in this paper should generalize to $SL(n, \mathbb Z)$ and perhaps to all congruence subgroups. It will require deeper studies on unital lattices in $\mathbb R^n$. Even though the precise statement for $\mc F_{\epsilon}$ will be difficult to write down, the Siegel sets of the same type may be obtained which will provide an equivalent norm for automorphic representations. Hence $\mc F_{\epsilon}$ or related Siegel set is potentially useful in the study of automorphic functions. The fundamental set $\mc F_{\epsilon}$ may also be studied from a topological viewpoint. It is not clear whether it can offer anything new. \\
\\
Let $\ceil{x}$ be the ceiling function, namely the smallest integer bigger than or equal to $x$.
\section{ Parametrization of $SL(2, \mathbb R)$: Setup}
Let $g \in SL(2, \mathbb R)$. Let  
$$g = (u,v), \qquad u =\arr{u_1 \\ u_2}, \qquad v= \arr{v_1 \\ v_2}.$$
Let
$$u^{\perp}=\arr{-\frac{u_1}{{u_1}^2+{u_2}^2} \\ \frac{u_2}{{u_1}^2+{u_2}^2}}.$$
Then $\|u ^{\perp} \|= \|u \|^{-1}$ and $v= u^{\perp}+ t u.$
The Iwasawa decomposition of $SL(2, \mathbb R)$ is given by 
$$g= (\frac{u}{\|u \|}, \frac{u^{\perp}}{\| u^{\perp} \|}) \arr{\|u \| & 0 \\ 0 & \| u^{\perp} \|} \arr{1 & t \\ 0 & 1}.$$
We may write this decomposition traditionally as $g=k a n$
with 
$$k=\arr{ \cos \theta & - \sin \theta \\ \sin \theta & \cos \theta}, \qquad 
a=\arr{\|u \| & 0 \\ 0 & \| u^{\perp} \|}, \qquad n=\arr{1 & t \\ 0 & 1}.$$
Here $ \theta \in \mathbb R/ 2 \pi \mathbb Z$ and  $t \in \mathbb R$. We may abuse the notation by writing $a= \|u \|$. Notice that $ \| u^{\perp} \|= 
\|u\|^{-1}=a^{-1}$ and
$$(\frac{u}{\|u \|}, \frac{u^{\perp}}{\| u^{\perp} \|})=\arr{\frac{u_1}{\|u \|} & -\frac{u_2}{\| u \|} \\ \frac{u_2}{\| u \|} & \frac{u_1}{\|u \|}} \in SO(2).$$
Since $v= u^{\perp}+ t u$, $\langle v, u \rangle= t \langle u, u \rangle$. Hence
$t = \frac{\langle v, u \rangle}{\langle u, u \rangle}$. \\
\\
Fix $K=SO(2)$ and $G=SL(2, \mathbb R)$.  We also have a variant of the Iwasawa decomposition $G= KNA$. The advantage of $KNA$ decomposition is that the product of the invariant measures on $K$, $N$ and $A$ is an invariant measure of $G$. To distinguish $KNA$ decomposition from $KAN$ decomposition, we write $g=k n(T) a$ in contrast with $g=k a n(t)$. In both decompositions, $k$ and $a$ remain the same. The parameter $T$ is related to $t$ by $T=a^2 t$. In fact, $T=\langle v, u \rangle$. \\
\\
In any case, all parameters pertaining to $g \in G$, $u, v, u^{\perp},t, a, T, \theta, k, n $ should be understood  as 
$$u(g), v(g), u^{\perp}(g),t(g), a(g), T(g), \theta(g), k(g), n(g)$$
for a fixed $g$.

\commentout{\subsection{Upper half plane $\mathbb H$}
The symmetric space $K \backslash G$ can be realized on the upper half plane $\mathbb H$}
\section{Parametrization of unital lattice in $\mathbb R^2$: $K$-invariant view}
Let $\Gamma=SL(2, \mathbb Z)$.
By a rank $2$ lattice in $\mathbb R^2$, we meant any additive subgroup of $\mathbb R^2$ with $2$ linearly independent generators. We call it a unital lattice if any two generators of the lattice span a parellelogram of area $1$. Denote the space of unital lattice in $\mathbb R^2$ by $\mc U$. Let $\mc L \in \mc U$. Then any pair of positively oriented generators $(u, v)$ of $\mc L$ corresponds an element $(u,v)$ of $G$ in a one-to-one fashion. Hence, the group $G$ parametrizes all positively oriented generators of all $\mc L \in \mc U$. Let the group $SL(2, \mathbb Z)$ acts on each lattice $\mc L \in \mc U$ by changing the generators:
$$(u, v) \rightarrow  (pu + qv, r u+ s v)=(u,v) \arr{p & r \\ q & s}, \qquad \forall \,\,\, \arr{p & r \\ q & s} \in \Gamma.$$
thus permuting the lattice points in $\mc L$.   Clearly any two positively oriented generators of $\mc L$ differs by an action of $\gamma \in \Gamma$ and vice versa.  Hence $\mc U$ can be identified with $G/\Gamma$ with $\Gamma$ acting from the right.  We have the natural projection 
$$\pi: G \rightarrow \mc U.$$
The fiber will be the matrices $\{ (u,v) \}$ given by any two positively oriented  generators $(u,v)$ of $\mc L $.
Furthermore, $\mc  U$ has a natural fibration $\tilde{\mc U}$:
$$ \mathbb Z^2 \rightarrow \tilde{\mc U }\rightarrow \mc U$$
with the group $\Gamma$ acting on the fiber as additive group automorphisms. \\
\\
Now we seek to parametrize $\mc U$. This is more or less equivalent to finding a fundamental set of $G/\Gamma$ in the classical sense. For a comprehensive account of fundamental set for reductiove groups, See \cite{borel} and the references therein. The fundamental sets in \cite{borel} are in fact the fundamental sets of $K \backslash G/\Gamma$. They can be viewed as a $K$-invariant set in  
$ G/\Gamma$ by the  pullback map of the projection $G \rightarrow K \backslash G$. Due to the action of the nontrivial center of $G$, the pullback will be a double cover of $G/\Gamma$ in the interior. In order to provide insight for the construction of non $K$-invariant fundamental set, we shall now review the basic ideas of the reduction theory in a $K$-invariant fashion. \\
\\
For each $\mc L \in \mc U$, let $\{ l^{(i)}  \}$ be the lattice points in $\mc L$. Define
$\Phi(\mc L)$ to be the minimal value of $\| l^{(i)} \|$ with $l^{(i)} \neq 0$. Fix 
a $u(\mc L) \in \mc L$ such that $\| u(\mc L) \|= \Phi(\mc L)$. Now $u(\mc L)$ is at least double-valued because $\| -u(\mc L)\|=\|u(\mc L)\|$. As we shall see later, there are  lattices that have more than two $u(\mc L)$. But these lattices will be of codimension one and with Haar measure zero. Hence it would not impact the general $L^p$ theory of $G/\Gamma$.  \\
\\
Now fix an $\mc L \in \mc U$. We shall write $u$ for $u(\mc L)$.
Because that $\| u \|$ is minimal among all lattice points in $\mc L$, there exists a $v$ such that
$$\mc L= \mathbb Z u + \mathbb Z v.$$
Since $\mc L$ is unital, $v= \pm u^{\perp}+t u$ and $\| v \|^2= \|u\|^{-2}+t^2 \| u \|^2$. Since $\|u \|$ is minimal, we have
$$\|u \|^{-2}+t^2 \|u \|^2 \geq \|u\|^2.$$
It follows that $t^2 \geq 1- \|u\|^{-4}$. \\
\\
Since the parameter $t \in \mathbb R/\mathbb Z$, we choose $t \in [-\frac{1}{2}, \frac{1}{2}]$.
Then $\frac{1}{4} \geq t^2 \geq 1-\|u\|^{-4}$. We obtain $\|u \| \leq (\frac{4}{3})^{\frac{1}{4}}$.
Now we define $\mc F$ to be the collection of $g=k a n$ with the properties
$$ a^2 \leq \frac{2}{\sqrt{3}}, \qquad t^2 \geq 1- a^{-4} , \qquad t^2 \geq \frac{1}{4}.$$
The following is well-known.
\begin{thm}
The natural map $\pi|_{\mc F}: \mc F \subset G \rightarrow G/\Gamma$ is surjective. It is two-to-one in the interior of $\mc F$ and finite on the boundary of $\mc F$.
\end{thm}
Proof: We already showed that if $\| u \|$ is minimal in $\mc L$, then 
$$ a^2 \leq \frac{2}{\sqrt{3}}, \qquad t^2 \geq 1- a^{-4} ,\qquad t^2 \geq \frac{1}{4}.$$
The converse is true, but not quite trivial. Suppose tha $g=kan$ with $a, t$ satisfy the above properties. For each $m \in \mathbb Z$, define the line
$$ u^{(m)}= m u^{\perp}+ \mathbb R u.$$
Then $\mc L \subseteq \sqcup_{m \in \mathbb Z} u^{(m)}$. If $|m| \geq 2$, then any $w \in u^{(m)}$ satisfies
$$\|w \|^2 \geq 4 \| u^{\perp} \|^2 = 4 \| u \|^{-2} = 4 a^{-2} \geq 4 \frac{\sqrt{3}}{2}= 2 \sqrt{3} > \frac{2}{\sqrt{3}} \geq \|u \|^2.$$
For $m= \pm 1$, $(u,v)$ or $(u, -v)$ is in $G$. Then $v= \pm u^{\perp}+ t u$. Our earlier discussion showed that $ \| v \| \geq \|u \|$ for any $v \in u^{(\pm 1)} \cap \mc L$. For $m=0$, $w \in u^{(0)} \cap \mc L$ means $w=k u$ for some $k \in \mathbb Z$. Hence $\| w \| \geq \|u \|$ unless $w=0$. \\
\\
Combining with the cases $m=0$ or $m=\pm 1$, we see that the minimal vector $u$ is unique up to a $\pm $ sign if $\frac{1}{4} > t^2 > 1- \| u\|^{-4}$.  Therefore  $\pi|_{\mc F}: \mc F \rightarrow G/\Gamma$ is two-to-one in the interior of $\mc F$.  Over the boundary of $\mc F$, it is not hard to see that $\deg(\pi|_{\mc F}) \leq 6$. When $t^2=1-a^{-4} $, the vectors $u$ and $v$ form an isosceles triangle and the degree of $\pi_{\mc F}$ over these points is $4$ with one exception: $t^2=1-a^{-4}= \frac{1}{4}$. This happens when $\mc L$ is generated by an equilateral triangle and $deg=6$ in this case. Of course, over $t=\pm \frac{1}{2}$, the degree of $\pi|_{\mc F}$ is also $4$ with the equilateral triangle case as the exception. $\Box$ 

\begin{cor}\label{one}
Let $\mc L \in \mc U$. If $\Phi(\mc L) < 1$, then there are only two vectors $\pm u \in \mc L$ such that $\| u \|=\Phi(\mc L)$.
\end{cor}
The boundary of $\mc F$ can be seen more easily on the upper half plane model of the symmetric space $K \backslash G$. The curve $t^2=1-a^{-4}$ corresponds to  a segment of the unit circle  and $t=\pm \frac{1}{2}$ corresponds to  two straight lines $x=\pm \frac{1}{2}$. Also if we use the projective group $PSL(2, \mb R)$, then map
$\pi|_{\mc F}$ will  be one-to-one in the interior of $\mc F$ and have degree at most $3$ over the boundary.\\
\\
Now we shall divide $\mc F$ into two regions:
$$\mc F^1= \{ g=kan \mid a \leq 1, |t| \leq \frac{1}{2} \};$$
$$\mc F^2= \{ g=kan \mid 1 \leq a \leq (\frac{4}{3})^{\frac{1}{4}},  1-a^{-4} \leq t^2 \leq \frac{1}{4} \}.$$
Then $\mc F =\mc F^1 \cup \mc F^2$.

\section{non $K$-invariant Parametrization of unital lattices in $\mathbb R^2$: main result}
When we parametrize $\mc U$ in the last section,  we allow the parameter $\theta \in \mathbb R/ 2 \pi \mathbb Z$ to be arbitrary. Fix $0 < \epsilon \leq \frac{\pi}{6}$. We now consider only $k$ with
$|\theta| \leq \epsilon$. Let
$C_{\epsilon}$ be defined as the open cone
$$\{ \arr{x \\y}: x >0, |\frac{y}{x}| < \tan \epsilon \} \subseteq \mathbb R^2. $$
Let $B(r)= \{ \arr{x \\ y} \in \mathbb R^2: x^2+y^2 \leq r^2 \}$. \\
\\
Fix $\mathcal L \in \mc U$. Let
$$\Phi_{\epsilon}(\mc L)= \min \{ \|u \|: u \in \mc L \cap C_{\epsilon} \}.$$
\begin{lem} $\Phi_{\epsilon}(\mc L)$ exists.
\end{lem}
Proof: Pick a generator $v \in \mc L$. Consider the lines $\{v^{(m)}: m \in \mathbb Z \}$. There are infinitely many $m_i$ with
$v^{(m_i)} \cap C_{\epsilon} \neq \emptyset$. If one of $v^{(m_i)} \cap C_{\epsilon}$is of infinite length, then there must be infinitely many lattice points on $v^{(m_i)} \cap C_{\epsilon}$. If not, the lengths of $v^{(m_i)} \cap C_{\epsilon}$, in ascending order, are of arithmetic progression, thus go to infinity. Recall that $\mc L \subseteq \sqcup v^{(m)}$ and $\mc L \cap v^{(m_i)}$ are equally spaced for all $m_i$. There must be infinitely many lattice points on these line segments. This show that the set $\{ \|u \|: u \in \mc L \cap C_{\epsilon} \}$ is infinite. Therefore  $I=\inf \{ \|u \|: u \in \mc L \cap C_{\epsilon} \}$ exists. \\
\\
To show that $\Phi_{\epsilon}(\mc L)$ exists, we must show that this infimium is a minimum. If not, there are infinitely many
lattice point $w_i$ such that $\| w_i \| \in [0, I+1]$. Then there must be infinitely many lattice points in the ball $B(I+1)$. This is not possible because $B(I+1)$ is compact. $\Box$.\\
\\
Now fix $ u \in \mc L$ such that $\| u \|=\Phi_{\epsilon}(\mc L)$. Strictly speaking, we should write $u_{\epsilon}(\mc L)$ for $u$. To simply our notation, we may write  $u$ or $u_{\epsilon}$, with the understanding  $\epsilon$ and $\mc L$ are fixed.
With $u_{\epsilon}$ selected for $\mc L \in \mc U$, we have the $k$ parameter  $(\cos \theta, \sin \theta)=u/\|u \|$ and $a = \|u \|$. For the $t$ parameter, we choose $t \in \mathbb R/\mathbb Z$.\\
\\
Define $$\mc F^1_{\epsilon}= \{ g=kan: |\theta| < \epsilon, a^2 \leq \csc(\epsilon+|\theta|), |t| \leq \frac{1}{2} \};$$
$$\mc F^2_{\epsilon}=\{ g=kan: |\theta| < \epsilon, \csc(\epsilon+|\theta|) \leq a^2 \leq \csc(\epsilon-|\theta|), t \in [-\sgn(\theta)a^{-2}\cot(\epsilon+|\theta|), \sgn(\theta)(1-\sqrt{1-a^{-4}})] \};$$
$$ \mc F^3_{\epsilon}=\{ g=kan: |\theta| < \epsilon, \csc(\epsilon-|\theta|) \leq a^2 \leq \cot(\epsilon-\theta)+\cot(\epsilon+\theta), t \in [1-a^{-2} \cot(\epsilon+\theta), a^{-2} \cot(\epsilon-\theta)] \}; $$
$$ \mc F^4_{\epsilon}=\{ g=kan: |\theta| < \epsilon, \csc(\epsilon-|\theta|) \leq a^2 , t \in [\sqrt{1-a^{-4}}-1, 1-\sqrt{1-a^{-4}}] \}.$$
Let $\mc F_{\epsilon}= \cup_{i=1}^4 \mc F^i_{\epsilon} \subseteq G$.
\begin{thm}  The natural map $\pi|_{\mc F_{\ep}}: \mc F_{\epsilon} \subseteq G \rightarrow G/\Gamma$ is surjective. $\pi|_{\mc F_{\ep}}$ is injective on the interior of $\mc F_{\ep}$ and finite on the boundary of $\mc F_{\ep}$.
\end{thm}
The condition that $\epsilon \leq \frac{\pi}{6}$ allows us to claim that $$\csc(\epsilon-|\theta|)\leq \cot(\epsilon-|\theta|)+\cot(\epsilon+|\theta|).$$ Hence $\mc F_{\ep}^3$ is not empty.
\section{Parametrization of $\mc U$:  necessary and sufficient condition for $t$}
Let $\mc L$ be a unital lattice in $\mathbb R^2$. Fix a $u$ in $\mc L \cap C_{\epsilon}$ such that
$\|u\|= \Phi_{\epsilon}(\mc L)$. $u$ must be primitive, namely there is no $u^{\prime} \in \mc L$ such that $u = q u^{\prime}$ with $|q| >1$. Then all lattice points of $\mc L$ must lay on one of the lines $u^{(m)}$ with $m \in \mathbb Z$. The lattice points on the line $u^{(m)}$ must be of the form $m u^{\perp}+(t+q) u$ for a   $t \in \mathbb R$ and $q \in \mathbb Z$. Hence $t$ is  in $\mathbb R/\mathbb Z$. In this section, we choose $t \in [0,1]$. Let $v=u^{\perp}+ t u$. Then $\mc L= \mathbb Z u+ \mathbb Z v$. Clearly, $u$ and $t$ parametrizes $\mc L$. 
For each $(a, \theta)$ determined by $u$, we need to find the range of $t$, i.e., a necessary and sufficient condition for $t$ such that \\
\\
{\bf for any $l \in \mc L \cap C_{\epsilon}$, $\| l\| \geq \|u \|$}. \\
\\
Let us call this property $\epsilon$. In contrast to the $K$-invariant parametrization where $\|u \|$ is bounded from above, 
$\| u_{\epsilon}(\mc L)\|$ is unbounded. \\
\\
 Define the (vertical) stripe
$$S= \{ s u^{\perp}+ t u \mid s \in \mathbb R, t \in (0,1) \}.$$
\begin{lem}\label{verti}
Let $v=u^{\perp}+ t u \in \mc L$ with $t \in [0,1]$. Then property $\epsilon$ is equivalent to the condition that
for any $l \in C_{\epsilon} \cap \mc L \cap S$, $\|l \| \geq \|u\|$.
\end{lem}
Proof: For $w= s u^{\perp}+ t u$ with $t \leq 0$, $w \notin C_{\epsilon}$.  For $w=s u^{\perp} + t u$ with $t \geq 1$, $\|w \| \geq \|u\|$. Hence we are only concerned with $w \in C_{\epsilon} \cap \mc L$ with $w=s u^{\perp}+t u \,\,\, (t \in (0,1))$. The equivalence with property $\epsilon$ follows immediately. $\Box$ \\
\\
Even though our lemme provide a necessary and sufficient condition, it is hard to manage the lattice points in
$ S \cap C_{\epsilon}^{0}$. If we fix a lattice point $v= u^{\perp}+ t u$, then the number of $ pu + qv $ in the triangle $S \cap C_{\epsilon}^{0}$ may depend on how $t$ is located near the rational points $\mathbb Q$. This turns out to be a difficult problem. \\
\\
Instead, we shall consider  a necessary condition, namely, for all lattice points $l \in u^{(\pm 1)} \cap \mc L \cap C_{\epsilon}$, $\|l \| \geq \|u \|$.
\begin{lem}\label{trange} Let $u_{\ep}(\mc L)= a\arr{\cos \theta \\ \sin \theta}$ with $\theta \in [0, \epsilon)$. Fix $t_{\epsilon}(\mc L) \in [0,1]$. Write $u, t$ for $u_{\ep}(\mc L),t_{\epsilon}(\mc L)$. 
\begin{enumerate}
\item For any $l \in \mc L \cap C_{\epsilon} \cap u^{(1)}$, $\| l \| \geq \|u \|$ if and only if one of the following is true
\begin{enumerate}
\item If $a^2 \leq \csc(\epsilon-\theta)$, then $t \in [0,1]$;
\item If $a^2 \geq \csc(\epsilon-\theta)$, then $ t \in [0, a^{-2} \cot(\epsilon-\theta)] \cup [\sqrt{1-a^{-4}}, 1]$.
\end{enumerate}
\item
 For any $l \in \mc L \cap C_{\epsilon} \cap u^{(-1)}$, $\| l \| \geq \|u \|$ if and only if one of the following is true
\begin{enumerate}
\item If $a^2 \leq \csc(\epsilon+\theta)$, then $t \in [0,1]$;
\item If $a^2 \geq \csc(\epsilon+\theta)$, then $ t \in [1- a^{-2} \cot(\epsilon+\theta), 1] \cup [0, 1-\sqrt{1-a^{-4}}]$.
\end{enumerate}
\end{enumerate}
\end{lem}
Proof: Let $v= u^{\perp}+ t u$. By the proof of Lemma \ref{verti}, we only need to find those $t \in [0,1]$ such that either $v \notin C_{\ep}$, or $ v \in C_{\ep}$ and $\|v \| \geq a$. By basic trigonometry, 
$$v \notin C_{\ep} \Longleftrightarrow \frac{\| t u \| }{\|u^{\perp} \|} \leq \cot(\ep-\theta) \Longleftrightarrow
t \leq a^{-2} \cot(\ep-\theta).$$
On the other hand, if $v \in C_{\ep}$,
$$\|v \|^2=t^2 \|u\|^2+\|u^{\perp}\|^2= a^2 t^2+a^{-2} \geq a^2= \|u \|^2 \Longleftrightarrow t^2 \geq 1-a^{-4}.$$
Hence either $t \in [0, a^{-2} \cot(\epsilon-\theta)] \cap [0,1]$ or $t \in [\sqrt{1-a^{-4}}, 1]$. If $a^{-2} \csc (\ep-\theta) \geq 1$, the union of these two sets is $[0,1]$. If $a^{-2} \csc (\ep-\theta) \leq 1$, then the union of these two sets is 
$[0, a^{-2} \cot(\epsilon-\theta)] \cup [\sqrt{1-a^{-4}}, 1]$. The first statement is proved.\\
\\
For the second statement, let $v=-u^{\perp}+(1-t)u$. Then $1-t$ should satisfy similar inequalities for the angle $\epsilon+\theta$. The second statement follows immediately. $\Box$\\
\\
Certainly, in order that $(u_{\ep}, t_{\ep})$ parametrizes $\mc L$, both conditions $(1)$ and $(2)$ must be met. We have

\begin{cor}Fix a $u \in \mc L$ such that $\|u\|=\Phi_{\epsilon}(\mc L)$. Suppose that $\theta \in [0, \ep)$ and $t_{\ep}(\mc L) \in [0,1]$. Then 
$$ t_{\epsilon}(\mc L)  \in ([0, a^{-2} \cot(\epsilon-\theta)] \cup [\sqrt{1-a^{-4}}, 1]) \cap ([1- a^{-2} \cot(\epsilon+\theta), 1] \cup [0, 1-\sqrt{1-a^{-4}}]) \cap [0,1].$$
\end{cor}
What is difficult and perhaps surprising is that these conditions turn out to be sufficient.\\
\\ 
We consider the horizontal direction (along $u$). Define $u^{+}$ to be the half plane $\{ w \in \mathbb R^2: \langle w, u^{\perp} \rangle > 0 \}$ and $u^{-}$ to be 
$\{ w \in \mathbb R^2: \langle w, u^{\perp} \rangle < 0 \}$. Then $\mathbb R^2= u^{+} \sqcup u^{(0)} \sqcup u^{-}$.
\begin{lem}\label{trange1}  Suppose that $u= a \arr{\cos \theta \\ \sin \theta}$ with $\theta \in [0, \epsilon)$ and $a >1$. If 
$$ t \in ([0, a^{-2} \cot(\epsilon-\theta)] \cup [\sqrt{1-a^{-4}}, 1]) \cap ([1- a^{-2} \cot(\epsilon+\theta), 1] \cup [0, 1-\sqrt{1-a^{-4}}]) \cap [0,1],$$
then the lattice $\mc L$ generated by $u$ and $u^{\perp}+t u$ is in $\mc U$ and satisfies the property that $\|u \|= \Phi_{\epsilon}(\mc L)$, the minimal norm of all lattice points in $\mc L \cap C_{\epsilon}$.
\end{lem}
{\bf Remark}: When $a \leq 1$, Lemma \ref{trange} says that $t \in \mathbb R/\mathbb Z$ can be arbitrary. This also turns out to be sufficient. More precisely,  the set
$$\{ \theta \in (-\ep, \ep), \qquad a \leq 1, \qquad t \in \mathbb R/\mathbb Z \}$$
is already contained in the fundamental set $\mc F$. The corresponding $u$ satisfies $\|u \|=\phi(\mc L)=\phi_{\ep}(\mc L)$. \\
\\
Proof: The key here is to treat $\mc L \cap u^{\pm}$ separately. Consider the half plane $u^+$. Set $$t \in [0, a^{-2} \cot(\epsilon-\theta)] \cup [\sqrt{1-a^{-4}}, 1]$$ and $v=u^{\perp}+ t u$. Define the line
$$v^{(m)}= m u + \mathbb R v, \qquad (m \in \mathbb Z).$$
The lattice $\mc L$ is contained in the union of the lines $v^{(m)}$ with $m \in \mathbb Z$. In fact, we have
$$\mc L \cap u^+= \mathbb Z^+ v+ \mathbb Z u.$$
\begin{enumerate}
\item Suppose that $ t \in [0, a^{-2} \cot(\epsilon-\theta)]$. Then $v =u^{\perp}+t u$ lies outside the cone $C_{\ep}$. For $m \leq 0$, $v^{(m)} \cap u^{+}$ lies entirely outside the cone $ C_{\ep}$. It suffice to show that for $m \in \mathbb Z^+$, $w \in v^{(m)} \cap u^+ \cap \mc L$, $ \| w \| \geq \|u \|$. Observe that $w = m u+j v$ for  $j \in \mathbb Z^+$ and  
$$\| m u + j v \|=\|(m+j t)u+ j u^{\perp} \| \geq \| u\|.$$ Hence for any $w \in C_{\epsilon} \cap u^{+} \cap \mc L$, $\| w \| \geq \|u \|$.

\item Suppose that $t^{\prime} \in [\sqrt{1-a^{-4}}, 1]$. Let $t=t^{\prime}-1$. Then $t \in [\sqrt{1-a^{-4}}-1, 0]$.  
Let $v=u^{\perp}+ t u $. Consider the semilattice $\mathbb Z^+ v + \mathbb Z u$. This is the same semilattice as if we use $v=u^{\perp}+ t^{\prime} u  $. 
\begin{enumerate}
\item for $m \leq 0$ the lattice points in $v^{(m)} \cap u^{+}$  are of the form $mu+ \mathbb Z^{+} v$. They lay entirely outside $C_{\ep}$. 
\item For $m > 1$ and $m, j \in \mathbb Z^+$, we have $m u+j v=(m+jt)u + j u^{\perp}$. If $|m+jt| \geq 1$, then
 $\| mu+j v \| \geq \|u\|$ automatically. If $ m+jt \in (-1, 1)$, then $j t < 1-m$. Hence 
 $$j > \frac{1-m}{t}=\frac{m-1}{-t} \geq \frac{1}{1-\sqrt{1-a^{-4}}} > a^4.$$
 It follows that $\| m u + j v \| \geq \|j u^{\perp} \| > a^3 > a =\| u \|$.
 \item For $m=1$, we would like to show that for $j \geq 1$,
  $$\| u+j v \|^2=\|(1+jt) u+ j u^{\perp}\|^2=(1+jt)^2a^2+j^2 a^{-2} \geq a^2=\|u\|^2.$$
 This is equivalent to 
 $(1+jt)^2 \geq 1-j^2 a^{-4}$.
 This follows from the fact that
 $$2t +j(t^2+a^{-4}) \geq 2t+ t^2+a^{-4}=(1+t)^2-1+a^{-4} \geq 0.$$
 \end{enumerate}
  Hence for any $w \in C_{\epsilon} \cap u^{+} \cap \mc L$, $\| w \| \geq \|u \|$.
  \end{enumerate}
 Next consider the half plane $u^-$. Set $$t^{\prime} \in [0, a^{-2} \cot(\epsilon+\theta)] \cup [ \sqrt{1-a^{-4}}, 1]$$ and $v^{\prime}=-u^{\perp}+ t^{\prime} u$.
 By essentially the same argument, 
 $$\|u \| \leq  \min \{\|w\|: w \in u^{-} \cap \mc L \cap C_{\ep} \}.$$ 
 The only difference is that the angle $\epsilon-\theta$ becomes $\epsilon+\theta$ because in $u^{-}$, the angle between the boundary of $C_{\ep}$ and $u$ is $\epsilon+\theta$. Notice this angle is less than $\pi/3$. Since $\det(u,v^{\prime})=-1$ in this setting, we switch back to the positive orientation  and obtain $v=u^{\perp}+ (1-t^{\prime}) u$. 
Since $t \in \mathbb R/\mathbb Z$ is chosen to be in $[0,1]$, we have
$$t \in  ([ 1-a^{-2} \cot(\epsilon+\theta),1] \cup [0, 1-\sqrt{1-a^{-4}}]) $$
if $a^{-2}\cot(\ep+\theta) \leq \sqrt{1-a^{-4}}$; $t \in [0,1]$ if $a^{-2}\cot(\ep+\theta) \geq \sqrt{1-a^{-4}}$. \\
\\
Observe that the range of $t$ we obtain is precisely the range of $t$ specified in Lemma \ref{trange1}. 
$\Box$ \\
\\

\section{ Proof of the main result }
We shall prove our main theorem for 
 $\theta \in [0, \ep)$. For $\theta \in (-\ep, 0]$, 
 the proof is similar.\\
 \\
 Let $\theta \in [0, \ep)$. If $a=\Phi_{\ep}(\mc L) \leq 1$, by Lemma \ref{trange}, $t \in [0,1]$. Conversely, for any such pair of $(a, t) \in (0,1] \times [0,1)$, there is one unique lattice with generators
$(u,v)$ satisfying $\|u \|=a$ and $v=u^{\perp}+t u$. In fact, this part of $\mc F_{\epsilon}$ overlaps with a small section of $K$-invariant $\mc F^1$. \\
\\
We shall now figure out precisely the range of $t$ for $a > 1$. Write
$$ I_{a, \theta}=([0, a^{-2} \cot(\epsilon-\theta)] \cup [\sqrt{1-a^{-4}}, 1]) \cap ([1- a^{-2} \cot(\epsilon+\theta), 1] \cup [0, 1-\sqrt{1-a^{-4}}]) \cap [0,1].$$
\begin{enumerate}
\item If  $a^2 \leq \csc(\epsilon+\theta)$, then $a^2 \leq \csc(\ep-\theta)$. We have
$$a^{-2} \cot(\ep-\theta) \geq a^{-2} \cot(\ep+\theta) \geq \sqrt{1-a^{-4}}.$$ Hence $t \in [0,1]$. Combined with the $a \leq 1$ case, $\mc L \in \mc F_{\epsilon}^1$. 
\item If  $\csc(\epsilon+\theta) \leq a^2 \leq \csc(\epsilon-\theta)$, then $$a^{-2} \cot(\ep-\theta) \geq   \sqrt{1-a^{-4}} \geq a^{-2} \cot(\ep+\theta).$$ We have $ t \in [1- a^{-2} \cot(\epsilon+\theta), 1] \cup [0, 1-\sqrt{1-a^{-4}}]$. Since $t \in \mathbb R/\mathbb Z$, we can make this set into a single interval: $t \in [-a^{-2}\cot(\epsilon+\theta), 1-\sqrt{1-a^{-4}}] $. We obtain $\mc L \in \mc F_{\epsilon}^2$. 
\item Suppose that $a^2 >\csc(\ep-\theta)$. Then 
$$ I_{a, \theta}=([0, a^{-2} \cot(\epsilon-\theta)] \sqcup [\sqrt{1-a^{-4}}, 1]) \cap ([1- a^{-2} \cot(\epsilon+\theta), 1] \sqcup [0, 1-\sqrt{1-a^{-4}}]).$$
We check that 
$$a^{-2} \cot(\epsilon-\theta) \geq 1-\sqrt{1-a^{-4}}, \qquad  \sqrt{1-a^{-4}} \geq 1- a^{-2} \cot(\epsilon+\theta)$$
by setting $a^{-2}=\sin \alpha$ and $0 < \alpha < \ep-\theta$. Then these two inequalities can be derived easily from the fact that
$$\cot(\ep-\theta) \geq \cot(\ep+\theta) \geq \cot( \frac{\pi}{3})=\frac{1}{\sqrt{3}} \geq \tan(\frac{\ep}{2}) \geq \tan(\frac{\alpha}{2}).$$
Hence 
$$I_{a, \theta}= [\sqrt{1-a^{-4}}, 1] \cup [0, 1-\sqrt{1-a^{-4}}] \cup ([0, a^{-2} \cot(\epsilon-\theta)] \cap [1- a^{-2} \cot(\epsilon+\theta), 1]).$$
If $a^2 \geq \cot(\ep+\theta)+\cot(\ep-\theta)$, then $I_{a, \theta}= [\sqrt{1-a^{-4}}, 1] \cup [0, 1-\sqrt{1-a^{-4}}]$. If $a^2 \leq \cot(\ep+\theta)+\cot(\ep-\theta)$, then
$$I_{a, \theta}= [\sqrt{1-a^{-4}}, 1] \cup [0, 1-\sqrt{1-a^{-4}}] \cup   [1- a^{-2} \cot(\epsilon+\theta),  a^{-2} \cot(\epsilon-\theta)].$$
Using a shift, we combine $[0, 1-\sqrt{1-a^{-4}}]$ with $[\sqrt{1-a^{-4}}-1, 0]$ and obtain $$ t \in [\sqrt{1-a^{-4}}-1, 1-\sqrt{1-a^{-4}} ].$$ Hence we obtain  $\mc F^3_{\ep}$ and  $\mc F^4_{\ep}$. We shall remark that $\cot(\ep+\theta)+\cot(\ep-\theta) \geq \csc(\ep-\theta)$ since $\ep \leq \frac{\pi}{6}.$ Therefore $\mc F_{\epsilon}^3 $ is nonempty.
\end{enumerate}
$\Box$ \\
\\
It is interesting to try to compare $\mc F_{\ep}$ with $\mc F$. Perhaps  the most distinctive feature is  that there is another "cusp" of the shape
$ |t| \leq 1-\sqrt{1-a^{-4}} \cong \frac{1}{2} a^{-4}$ as $a \rightarrow \infty$ in $\mc F_{\epsilon}$. \\
\\
We can now rewrite $\mc F_{\epsilon}$ in the $KNA$ decomposition.
$$\mc F^1_{\epsilon}= \{ g=kan: |\theta| < \epsilon, a^2 \leq \csc(\epsilon+|\theta|), |T| \leq \frac{1}{2} a^2 \};$$
$$\mc F^2_{\epsilon}=\{ g=kan: |\theta| < \epsilon, \csc(\epsilon+|\theta|) \leq a^2 \leq \csc(\epsilon-|\theta|), \sgn(\theta) T \in [-\cot(\epsilon+|\theta|), (a^2-\sqrt{a^4-1})] \};$$
$$ \mc F^3_{\epsilon}=\{ g=kan: |\theta| < \epsilon, \csc(\epsilon-|\theta|) \leq a^2 \leq \cot(\epsilon-\theta)+\cot(\epsilon+\theta), T \in [a^2- \cot(\epsilon+\theta),  \cot(\epsilon-\theta)] \}; $$
$$ \mc F^4_{\epsilon}=\{ g=kan: |\theta| < \epsilon, \csc(\epsilon-|\theta|) \leq a^2 , T \in [\sqrt{a^4-1}-a^2, a^2-\sqrt{a^{4}-1}] \}.$$
For $\mc F_{\ep}^3$ and $\theta \in [0,\epsilon)$, we can  shift the parameter $T$ to $(-\cot(\ep+\theta), \cot(\ep-\theta)-a^2)$.

\begin{thm}\label{big} Let $\ep \in (0, \frac{\pi}{6})$.
The fundamental set can be chosen inside
$$\{ \theta \in (-\ep, \ep), a \in (0, \infty), |T| \leq \cot(\ep+|\theta|)\} \subseteq G.$$
\end{thm}
Proof: Since $\ep+|\theta| <\frac{\pi}{3}$, in $\mc F_{\ep}^1$,
$$|T| \leq \frac{1}{2} a^2 \leq \frac{1}{2} \csc(\ep+|\theta|) < \cot(\ep+|\theta|).$$
In $\mc F_{\ep}^2$, since $a^2 \geq \csc(\ep+|\theta|)$,
$$a^2-\sqrt{a^4-1} \leq \cot(\ep+|\theta|).$$
Hence $|T| \leq \cot(\ep+|\theta|)$. 
In $\mc F^{3}_{\ep}$, if $\theta \in (-\ep, 0)$, then $T \in [a^2- \cot(\epsilon+\theta),  \cot(\epsilon-\theta)]$.  Clearly $|T| \leq \cot(\ep+|\theta|)$. If $\theta \in [0, \ep)$, then $T \in (-\cot(\ep+\theta), \cot(\ep-\theta)-a^2)$. We also have $|T| \leq \cot(\ep+|\theta|)$. In $\mc F_{\ep}^4$, we clearly have
$$a^2-\sqrt{a^4-1} \leq \cot(\ep+|\theta|).$$ 
Hence
  $\mc F_{\ep}$ can be put entirely inside   $$\{ \theta \in (-\ep, \ep), a \in (0, \infty), |T| \leq \cot(\ep+|\theta|) \}.$$
  When $a^2 \in [\csc(\ep+|\theta|), \cot(\ep+\theta)+\cot(\ep-\theta)]$, this is the best approximation. For $a$ small or big, we have two cusps: one  at $0$ and one at $\infty$. $\Box$ \\
\\
Notice that $\tan(\ep+|\theta|) \geq \tan (\ep)\geq \ep$. Hence $\cot(\ep+\theta) \leq \frac{1}{\ep}$.
In $\mc F_{\ep}$, there is  a "duality" between the $k$ parameter and $T$ parameter, namely as the size of $|\theta|$ shrinks, the range of $T$ can increase to size $\frac{1}{\epsilon}$. \\
\\
Now we shall give an integral inequality on functions on $G/\Gamma$. 
Since $\cot(\ep+|\theta|) \leq \cot(\ep) < \frac{1}{\ep}$, the following  theorem is  an immediate consequence of Theorem \ref{big}. 
\begin{thm} Let $f \in L^2(G/\Gamma)$. Then
$$\int_{-\ep}^{\ep} \int_{-\frac{1}{\ep}}^{\frac{1}{\ep}} \int_{0}^{\infty} |f(k_{\theta} n_{T} a)|^2 \frac{d a}{a} d T d \theta
\geq \| f\|_{L^2}^2.$$
Here $k_{\theta} n_{T} a$ is the $KNA$ decomposition and $\frac{d a}{a} d T d \theta$ is the $G$-invariant measure.
\end{thm}


\begin{thebibliography}{99}


\bibitem{ba} H. Bateman,{\it Higher Transcendental Functions  Vol.I-III,} McGraw-Hill Book Company, 1953.
\bibitem{borel} A. Borel {\it Automorphic forms on $SL(2)$}  Cambridge Tracts in Mathematics, 130. Cambridge University Press, Cambridge, 1997. 
\bibitem{borel19} A. Borel {\it Introduction to Arithmetic Groups}, American Mathematical Society, Providence 2019.
\bibitem{bh} A. Borel and Harish-Chandra, \lq\lq Arithmetic subgroups of Algebraic groups, \rq\rq {\it Annals of Math.} Vol. 75, (485-535) 1962.
\bibitem{go} D. Goldfeld {\it Automorphic Forms and L-functions for the Group $GL(n, \mathbb R)$}, Cambridge University Press, Cambridge 2006.
\bibitem{hc} Harish-Chandra {\it Automorphic Forms on Semisimple Lie Groups}, Notes by J. G. M. Mars, LNM 62, Springer-Verlag, 1968.
\bibitem{he1} H. He, \lq\lq  Representations of $ax+b$ group and Dirichlet Series,\rq\rq   J. Ramanujan Math. Soc., Vol. 36 2021 (73-84).
\bibitem{he2} H. He \lq\lq Certain $L^2$-norms on authomorphic representations of $SL(2, \mathbb R)$,
\rq\rq submitted.
\bibitem{knapp} A. Knapp {\it Representation theory of semisimple groups}, Princeton University Press 2002.
\bibitem{val} N.  Vilenkin {\it Special Functions and the Theory of Group Representations}, American Mathematical Society, Providence, 1983.
\end{thebibliography}
\end{document}